\documentclass{amsart}
\usepackage{ifthen, citesort}
\usepackage{amssymb}
\usepackage{epic}

\DeclareMathOperator{\tr}{tr}
\DeclareMathOperator{\divergence}{div}
\DeclareMathOperator{\divergenz}{div}

\DeclareMathOperator{\diag}{diag}
\DeclareMathOperator{\osc}{osc}
\def\R{\mathbb{R}}
\def\N{\mathbb{N}}
\def\S{\mathbb{S}}
\def\Z{\mathbb{Z}}
\def\dt{\frac{d}{dt}}
\def\fracp#1#2{\frac{\partial #1}{\partial #2}}

\def\theta{\vartheta}
\def\phi{\varphi}
\def\epsilon{\varepsilon}
\newcommand{\A}[1]{\ifthenelse{#1 = 2}{\lvert A\rvert^{#1}}{\tr A^{#1}}}
\newcommand{\Akl}[2]{\ifthenelse{#1 = 2}%
{\left(\lvert A\rvert^{#1}\right)^{#2}}%
{\left(\tr A^{#1}\right)^{#2}}}
\def\Az#1{\left(\A2\right)_{;\,#1}}
\def\B#1{\tr A^{#1}}
\def\hij{h_{ij;\,k}^2}
\def\hii{h_{ii;\,k}h_{jj;\,k}}
\def\li{\lambda_i}
\def\lj{\lambda_j}
\def\lk{\lambda_k}
\def\lx{\lambda_1}
\def\ly{\lambda_2}
\def\opr#1{\dt{#1}-F^{ij}{\left(#1\right)}_{;\,ij}}
\def\oprokl#1{\dt{#1}-F^{ij}{#1}_{;\,ij}}
\def\fracd#1#2{\displaystyle\frac{\displaystyle #1}{\displaystyle #2}}

\long\def\umbruch{{\displaybreak[1]}}

\newtheorem{theorem}{Theorem}[section]
\newtheorem{lemma}[theorem]{Lemma}
\newtheorem{proposition}[theorem]{Proposition}
\newtheorem{corollary}[theorem]{Corollary}

\theoremstyle{definition}

\theoremstyle{remark}
\newtheorem{remark}[theorem]{Remark}

\numberwithin{equation}{section}
\newcommand{\abs}[1]{\left\lvert#1\right\rvert}

\begin{document}

\title{Surfaces Expanding by the Inverse Gau{\ss} Curvature Flow}

\author{Oliver C. Schn\"urer}
\address{FU Berlin, Arnimallee 2-6, 14195 Berlin, Germany}
\curraddr{}
\email{Oliver.Schnuerer@math.fu-berlin.de}
\thanks{The author is a member of SFB 647/B3 ``Raum -- Zeit -- Materie''.}


\date{December 2004, revised February 2006.}

\dedicatory{}


\begin{abstract}
We show that strictly convex surfaces expanding
by the inverse Gau{\ss} curvature flow converge
to infinity in finite time. After appropriate rescaling, 
they converge to spheres. 
We describe the algorithm to find our main test function. 
\end{abstract}

\maketitle

\tableofcontents

\section{Introduction}

We consider a family of closed strictly convex surfaces $M_t$ in $\R^3$ 
that expand by the inverse Gau{\ss} curvature flow
\begin{equation}\label{flow eqn}
\dt X=\frac1K\nu.
\end{equation}
This is a parabolic flow equation. We obtain a solution
on a maximal time interval $[0,\,T)$, $0<T<\infty$. 
For $t\uparrow T$, the surfaces converge to infinity. After
appropriate rescaling, they converge to a round sphere.
We say that the surfaces $M_t$ converge to ``round spheres at infinity''.
The key step in the proof, Theorem \ref{K invers mon thm}, is to show that
\begin{equation}\label{mon groe}
\max\limits_{M_t}\left(\frac{(\lx-\ly)^2}{\lx^2\ly^2}\right)
\end{equation}
is non-increasing in time. \par
Here, we used standard notation 
as explained in Section \ref{nota sec}.

Our main theorem is
\begin{theorem}\label{main thm}
For any smooth closed strictly convex surface $M$ in $\R^3$, 
there exists a smooth family of surfaces $M_t$, $t\in[0,T)$,
solving \eqref{flow eqn} with $M_0=M$. For $t\uparrow T$,
the surfaces $M_t$ converge to infinity. The rescaled surfaces
$M_t\cdot(T-t)$ converge smoothly to the
unit sphere $\S^2$.
\end{theorem}

We will also consider other normal velocities for which
similar results hold. Therefore, we have to find 
quantities like \eqref{mon groe} that are monotone
during the flow and vanish precisely for spheres.
In general, this is a complicated issue. 
In order to find these test quantities, we used an algorithm
that checks, based on randomized tests, whether possible
candidates fulfill certain inequalities. These inequalities
guarantee especially that we can apply the maximum principle
to prove monotonicity. We used that algorithm only to propose
useful quantities. The presented proof does not depend on it. 
So far, all candidates proposed by the corresponding program 
turned out to be
appropriate for proving convergence to round spheres at infinity. 
In Table \ref{table}, we have collected some normal
velocities $F$ ($1^{\text{st}}$ column) and quantities
$w$ ($2^{\text{nd}}$ column) such that $\max_{M_t}w$ is
non-increasing in time for surfaces expanding with normal velocity $F$
$$\dt X=-F\nu.$$ It is common to use positive
functions $F$ for contracting surfaces. Thus the negative
sign corresponds to the fact that these surfaces expand.
In each case, we obtain
convergence to round spheres at infinity for smooth closed strictly 
convex initial surfaces $M_0$. 

\begin{table}
\def\platz{\raisebox{0em}[2.2em][1.5em]{\rule{0em}{2em}}}
$$
\begin{array}{|c||c|}\hline
-\fracd1K & \platz\fracd{(\lx-\ly)^2}{\lx^2\ly^2}\\\hline
-\fracd{H^2}{K^2} & \platz\fracd{(\lx-\ly)^2}{(\lx+\ly)\lx\ly}\\\hline
-\fracd{\A2}{K^2} & \platz\fracd{(\lx-\ly)^2}{(\lx+\ly)\lx\ly}\\\hline
-\fracd{H^3}{K^3} & \platz\fracd{(\lx+\ly)^6(\lx-\ly)^2}
{\left(\lx^2+\ly^2\right)\left(\lx^2+\lx\ly+\ly^2\right)\lx^3\ly^3}\\
\hline
\end{array}$$
\caption{Monotone quantities}\label{table}
\end{table}

There are many papers concerning convex hypersurfaces
contracting to ``round points'', i.\,e.\ the surfaces converge to 
a point, and, after appropriate rescaling, to a sphere, especially for
normal velocities homogeneous of degree one in the
principal curvatures, see e.\,g.\ \cite{HuiskenRoundSphere}
for motion by mean curvature. 
For normal velocities of higher homogeneity, strictly convex
hypersurfaces converge to a point \cite{TsoPoint}
and to a round point, if they are appropriately
pinched initially, e.\,g.\ \cite{Andrews2dnonconcave} 
and \cite{OSFelixH2} for surfaces. 
Convex surfaces contracting with normal velocities 
homogeneous of degree larger than one in the principal 
curvatures converge to round points without initial
pinching assumption, see \cite{AndrewsStones} for
the Gau{\ss} curvature flow and \cite{OSA2} for other
flow equations and strictly convex surfaces. 

Expanding flows of homogeneity minus one, i.\,e.\
flows of the form $\dt X=-F\nu$ with $F$ positive
homogeneous of degree minus one, were studied by
Claus Gerhardt and John Urbas
\cite{CGFlowSpheres,UrbasExpandJDG,UrbasExpandMZ}.
They obtain convergence to round spheres at infinity.
These results extend to negative homogeneities larger than
minus one. Note that solutions exist for $t\in[0,\,\infty)$.
There is a representation formula for solutions to the
inverse harmonic mean curvature flow of Knut Smoczyk 
\cite{KnutHarmonicMCF}. Gerhard Huisken and Tom Ilmanen
used the inverse mean curvature flow to prove the 
Penrose inequality \cite{HuiskenIlmanenPenrose}. This
was extended in \cite{HuiskenIlmanenIMCF2,%
SmoczykAsian2000}. 
There are also inhomogeneous flows for which solutions
converge to round spheres as $t\uparrow\infty$ 
\cite{MinkowskiFlow,ChowAsian,ChowNonhomogGauss,%
IvochkinaNehringTomi}.

Our paper concerns the expansion of surfaces by the inverse 
Gau{\ss} curvature flow. This flow is homogeneous of degree 
minus two in the principal curvatures. Solutions tend to
round points at infinity in finite time. 
We want to stress that we don't
have to assume any pinching condition for the initial surface.

The rest of this paper is organized as follows. In Section 
\ref{nota sec}, we explain our notation. 
Section \ref{mon sec} concerns the key step,
Theorem \ref{K invers mon thm}, the proof of the monotonicity of 
our test function during the flow.
We prove convergence to infinity
for appropriate points on the surfaces
in Section \ref{inf sec}. In Section \ref{conv sphere sec}, 
we obtain that the surfaces converge to infinity and, after
appropriate rescaling, to a round sphere in Hausdorff distance.
We improve this in Section \ref{smooth conv sec} and get smooth
convergence to a round sphere after rescaling. This finishes
the proof of our main theorem. A further improvement of this 
result is contained in Section \ref{impr conv sec}. There, we 
show that our surfaces converge in Hausdorff distance to a
family of spheres expanding by inverse Gau{\ss} curvature flow.

In Sections \ref{algor sec} and \ref{comp alg sec}, we describe
the computational aspects of our flow equation. 
We describe the algorithm that we used to find our monotone
quantities. In this paper, we compute complicated evolution
equations. We describe, how this can be done with a computer
algebra program. 

We also get convergence to round spheres at infinity for other
flow equations. These are considered in Section \ref{o nor vel}.
Finally, we derive the optimal expected convergence rate in Section
\ref{conv rate}. 

The author wants to thank Klaus Ecker at the Free University Berlin
for discussions and support, especially for telling us about
Aleksandrov reflection for parabolic equations. We also want to
thank J\"org H\"arterich, Gerhard Huisken, and Felix Otto for 
discussions concerning the optimal convergence rate.

\section{Notation}\label{nota sec}
We use $X=X(x,\,t)$ to denote the embedding vector of a manifold 
$M$ into $\R^3$ and $\dt X=\dot X$ for its total time derivative. 
Set $M_t:=X(M,\,t)\subset\R^3$.
We choose $\nu$ to be the outer unit normal vector of $M_t$. 
The embedding induces a metric $(g_{ij})$ and
a second fundamental form $(h_{ij})$. We use the Einstein summation
convention. Indices are raised and lowered with respect to the metric
or its inverse $\left(g^{ij}\right)$. The inverse of the second fundamental
form is denoted by $\left(\tilde h^{ij}\right)$. The principal
curvatures $\lx,\,\ly$ are the eigenvalues of the second fundamental
form with respect to the metric. A surface is called strictly convex,
if all principal curvatures are strictly positive.
We will assume this throughout the paper.

Symmetric functions of the principal
curvatures are well-defined, we will use the mean curvature
$H=\lx+\ly$, the square of the norm of the second fundamental form
$\A2=\lx^2+\ly^2$, $\B k=\lx^k+\ly^k$, and the Gau{\ss} curvature
$K=\lx\ly$. We write indices, preceded by semi-colons, e.\,g.\ $h_{ij;\,k}$, 
to indicate covariant differentiation with respect to the induced
metric. It is often convenient to choose coordinate systems such
that, at a fixed point, 
the metric tensor equals the Kronecker delta, $g_{ij}=\delta_{ij}$,
and $(h_{ij})$ is diagonal, $(h_{ij})=\diag(\lx,\,\ly)$, e.\,g.
$$\sum\lk\hij=\sum\limits_{i,\,j,\,k=1}^2\lk\hij
=h^{kl}h^i_{j;\,k}h^j_{i;\,l}
=h_{rs}h_{ij;\,k}h_{ab;\,l}g^{ia}g^{jb}g^{rk}g^{sl}.$$
Whenever we use this notation, we will also assume that we have 
fixed such a coordinate system. We will only use Euclidean coordinate
systems for $\R^3$ so that $h_{ij;\,k}$ is symmetric according to
the Codazzi equations.

A normal velocity $F$ can be considered as a function of $(\lx,\,\ly)$
or $(h_{ij},\,g_{ij})$. We set $F^{ij}=\fracp{F}{h_{ij}}$,
$F^{ij,\,kl}=\fracp{^2F}{h_{ij}\partial h_{kl}}$. 
Note that in coordinate
systems with diagonal $h_{ij}$ and $g_{ij}=\delta_{ij}$ as mentioned
above, $F^{ij}$ is diagonal. 
For $F=-K^{-1}$, we have $F^{ij}=K^{-1}\tilde h^{ij}=
K^{-1}\li^{-1}g^{ij}$. 

Recall, see e.\,g.\ \cite{HuiskenRoundSphere,OSKnutCrelle,OSMM,%
HuiskenPolden}, 
that for a hypersurface moving according to $\dt X=-F\nu$, we have
\begin{align}
\label{g evol}
\dt g_{ij}=&-2Fh_{ij},\displaybreak[1]\\
\label{h evol}
\dt h_{ij}=&F_{;\,ij}-Fh_i^kh_{kj},\displaybreak[1]\\
\dt\nu^\alpha=&g^{ij}F_{;\,i}X^\alpha_{;\,j},
\end{align}
where Greek indices refer to components in the ambient space $\R^3$.
In order to compute evolution equations, we will need the
Gau{\ss} equation and the Ricci identity for the second fundamental
form
\begin{align}
\label{Riem}
R_{ijkl}=&h_{ik}h_{jl}-h_{il}h_{jk},\displaybreak[1]\\
\label{Ricci}
h_{ik;\,lj}=&h_{ik;\,jl}+h^a_kR_{ailj}+h^a_iR_{aklj}.
\end{align}
We will also employ the Gau{\ss} formula and the Weingarten
equation
$$X^\alpha_{;\,ij}=-h_{ij}\nu^\alpha\qquad\text{and}\qquad
\nu^\alpha_{;\,i}=h^k_iX^\alpha_{;\,k}.$$

For tensors $A$ and $B$, $A_{ij}\ge B_{ij}$ means that 
$(A_{ij}-B_{ij})$ is positive definite.
Finally, we use $c$ to denote universal, estimated constants. 

\section{A Monotone Quantity}\label{mon sec}

\begin{theorem}\label{K invers mon thm}
For a family of smooth closed strictly convex surfaces $M_t$ 
in $\R^3$ flowing according to $\dot X=\frac1K\nu$, 
\begin{equation}\label{K invers quantity}
\max\limits_{M_t}\frac{(\lambda_1-\lambda_2)^2}
{4\lambda_1^2\lambda_2^2}
=\max\limits_{M_t}\frac{2\A2-H^2}{\left(H^2-\A2\right)^2}
\equiv\max\limits_{M_t}w
\end{equation}
is non-increasing in time.
\end{theorem}

An immediate consequence of this theorem is
\begin{corollary}
The only homothetically expanding smooth closed strictly convex
surfaces $M_t$, solving the flow equation $\dot X=\frac1K\nu$ in 
$\R^3$, are spheres.
\end{corollary}
\begin{proof}
The quantity $\frac{(\lx-\ly)^2}{\lx^2\ly^2}$ 
is negative homogeneous of degree two
in the principal curvatures and non-negative. If $M$ is 
homothetically expanding, Theorem \ref{K invers mon thm} implies that
$(\lx-\ly)^2=0$ everywhere. Thus $M_t$ is umbilic and 
\cite[Lemma 7.1]{Spivak4} implies that $M_t$ is a sphere.
\end{proof}

\begin{proof}[Proof of Theorem \ref{K invers mon thm}]
First and second derivatives of $F=-1/K$ with respect 
to the second fundamental form are given by
\begin{align}\label{F abl}
\begin{split}
F^{ij}=&\frac1K\tilde h^{ij},\\
F^{ij,\,kl}=&-\frac1K\left(\tilde h^{ij}\tilde h^{kl}+
\tilde h^{ik}\tilde h^{jl}\right).
\end{split}
\end{align}
We combine \eqref{g evol}, \eqref{h evol}, \eqref{Riem},
and \eqref{Ricci} in order to get the 
following evolution equations
\begin{align}\label{H evol}
\oprokl H=&\frac1K\left(H^2-3\A2\right)
-\frac1K\sum\frac1{\li\lj}\left(\hii+\hij\right)
\intertext{and}
\begin{split}\label{A2 evol}
\opr{\A2}=&\frac1K\left(2H\A2-6\A3\right)-\frac1K\sum\frac2\lk\hij\\
&-\frac1K\sum\frac{2\lk}{\li\lj}
\left(\hii+\hij\right).
\end{split}
\end{align}
For the rest of the proof, we consider a critical point of 
$\left.w\right|_{M_t}$ for some $t>0$, where $w>0$
and $\lx\neq\ly$.
It suffices to show that $w$ is non-increasing
at such a point. Then our theorem follows. 

We differentiate $w$
\begin{align*}
w=&\frac{-H^2+2\A2}{\left(H^2-\A2\right)^2},\umbruch\\
w_{;\,i}=&\frac{2H\left(H^2-3\A2\right)H_{;\,i}+2\A2\Az i}
{\left(H^2-\A2\right)^3},\umbruch\\
w_{;\,ij}=&\frac1{\left(H^2-\A2\right)^4}\left\{
2H\left(H^2-\A2\right)\left(H^2-3\A2\right)H_{;\,ij}\right.\\
&\quad\left.+2\A2\left(H^2-\A2\right)\Az{ij}\right.\\
&\quad\left.\left(6\left(H^2-\A2\right)^2-12H^2\left(H^2-3\A2\right)
\right)H_{;\,i}H_{;\,j}\right.\\
&\quad\left.+2\left(H^2+2\A2\right)\Az i\Az j\right.\\
&\quad\left.-12H\A2\left(H_{;\,i}\Az j+H_{;\,j}\Az i\right)\right\}.
\end{align*}
In a critical point of $w$ with $\lx\neq\ly$, we get
\begin{align}
0=&4\ly^2(\lx-\ly)h_{11;\,k}+4\lx^2(\ly-\lx)h_{22;\,k}
\quad\text{for }k=1,\,2,\umbruch\nonumber\\
h_{22;\,1}=&\frac{\ly^2}{\lx^2}h_{11;\,1}
\equiv a_1h_{11;\,1}\quad\text{and}\quad
h_{11;\,2}=\frac{\lx^2}{\ly^2}h_{22;\,2}
\equiv a_2h_{22;\,2}.\label{extremal cond}
\end{align}
The evolution equation of $w$ is
\begin{align*}
\left(H^2-\A2\right)^4\cdot&\left(\oprokl w\right)=\\
=&2H\left(H^2-\A2\right)\left(H^2-3\A2\right)\left(\oprokl H\right)
\umbruch\\
&+2\A2\left(H^2-\A2\right)\left(\opr{\A2}\right)\umbruch\\
&+\left(6\left(H^2-\A2\right)^2-12H^2\left(H^2-3\A2\right)\right)
\left(-F^{ij}H_{;\,i}H_{;\,j}\right)\umbruch\\
&+2\left(H^2+2\A2\right)\left(-F^{ij}\Az i\Az j\right)\umbruch\\
&-12H\A2\left(-F^{ij}\left(H_{;\,i}\Az j+H_{;\,j}\Az i\right)\right).
\end{align*}
We consider the terms on the right-hand side individually and use
\eqref{H evol}, \eqref{A2 evol}, and \eqref{extremal cond}
\begin{align*}
\oprokl H=&-\frac{2\left(\lx^2-\lx\ly+\ly^2\right)}{\lx\ly}\umbruch\\
&-\frac{2\left(\lx^2a_1^2+\lx\ly a_1^2+\lx\ly a_1+\ly^2\right)}
{\lx^3\ly^3}h_{11;\,1}^2\umbruch\\
&-\frac{2\left(\ly^2a_2^2+\ly\lx a_2^2+\ly\lx a_2+\lx^2\right)}
{\ly^3\lx^3}h_{22;\,2}^2\umbruch\\
=&-\frac{2\left(\lx^2-\lx\ly+\ly^2\right)}{\lx\ly}\umbruch\\
&-\frac{2(\lx+\ly)\left(\lx^2+\ly^2\right)}{\lx\ly}
\left(\frac1{\lx^5}h_{11;\,1}^2+\frac1{\ly^5}h_{22;\,2}^2\right),
\umbruch\\
2H\left(H^2-\A2\right)&\left(H^2-3\A2\right)\cdot
\left(\oprokl H\right)=\umbruch\\
=&-8(\lx+\ly)\left(\lx^2-\lx\ly+\ly^2\right)\lx\ly
\left(\oprokl H\right)\umbruch\\
=&16(\lx+\ly)\left(\lx^2-\lx\ly+\ly^2\right)^2\umbruch\\
&+16(\lx+\ly)^2\left(\lx^2-\lx\ly+\ly^2\right)\left(\lx^2+\ly^2\right)
\cdot\\
&\quad\cdot
\left(\frac1{\lx^5}h_{11;\,1}^2+\frac1{\ly^5}h_{22;\,2}^2\right)\umbruch\\
=&8(\lx+\ly)
\left(2\lx^4-4\lx^3\ly+6\lx^2\ly^2-4\lx\ly^3+2\ly^4\right)\umbruch\\
&+\left(4\lx^7\ly+4\lx^6\ly^2+4\lx^5\ly^3+8\lx^4\ly^4\right.\\
&\qquad\left.+4\lx^3\ly^5+4\lx^2\ly^6+4\lx\ly^7\right)\cdot\\
&\quad\cdot\frac4{\lx\ly}\left(\frac1{\lx^5}h_{11;\,1}^2
+\frac1{\ly^5}h_{22;\,2}^2\right),\umbruch\\
\opr{\A2}=&-\frac{2(\lx+\ly)\left(2\lx^2-3\lx\ly+2\ly^2\right)}{\lx\ly}
\umbruch\\
&-\frac{2\left(2\lx^2 a_1^2+2\lx\ly a_1^2+3\ly^2a_1^2
+2\lx\ly a_1+3\ly^2\right)}{\lx^2\ly^3}h_{11;\,1}^2\umbruch\\
&-\frac{2\left(2\ly^2a_2^2+2\ly\lx a_2^2+3\lx^2a_2^2
+2\ly\lx a_2+3\lx^2\right)}{\ly^2\lx^3}h_{22;\,2}^2\umbruch\\
=&-\frac{2(\lx+\ly)\left(2\lx^2-3\lx\ly+2\ly^2\right)}{\lx\ly}\umbruch\\
&-\frac{2\left(3\lx^4+2\lx^3\ly+2\lx^2\ly^2+2\lx\ly^3+3\ly^4\right)}
{\lx\ly}\cdot\\
&\quad\cdot\left(\frac1{\lx^5}h_{11;\,1}^2+\frac1{\ly^5}h_{22;\,2}^2\right),
\umbruch\\
2\A2\left(H^2-\A2\right)\cdot&\left(\opr{\A2}\right)=\umbruch\\
=&4\left(\lx^2+\ly^2\right)\lx\ly\left(\opr{\A2}\right)\umbruch\\
=&-8(\lx+\ly)\left(\lx^2+\ly^2\right)\left(2\lx^2-3\lx\ly+2\ly^2\right)
\umbruch\\
&-8\left(\lx^2+\ly^2\right)
\left(3\lx^4+2\lx^3\ly+2\lx^2\ly^2+2\lx\ly^3+3\ly^4\right)\cdot\umbruch\\
&\quad\cdot\left(\frac1{\lx^5}h_{11;\,1}^2+\frac1{\ly^5}h_{22;\,2}^2\right)
\umbruch\\
=&-8(\lx+\ly)\left(2\lx^4-3\lx^3\ly+4\lx^2\ly^2-3\lx\ly^3+2\ly^4\right)
\umbruch\\
&-\left(6\lx^7\ly+4\lx^6\ly^2+10\lx^5\ly^3+8\lx^4\ly^4+10\lx^3\ly^5
\right.\\
&\qquad\left.+4\lx^2\ly^6+6\lx\ly^7\right)\frac4{\lx\ly}
\left(\frac1{\lx^5}h_{11;\,1}^2+\frac1{\ly^5}h_{22;\,2}^2\right).
\end{align*}
For the remaining terms, we employ \eqref{F abl}
\begin{align*}
-F^{ij}H_{;\,i}H_{;\,j}=&-\frac{(a_1+1)^2}{\lx^2\ly}h_{11;\,1}^2
-\frac{(a_2+1)^2}{\ly^2\lx}h_{22;\,2}^2\umbruch\\
=&-\frac{\left(\lx^2+\ly^2\right)^2}{\lx\ly}
\left(\frac1{\lx^5}h_{11;\,1}^2+\frac1{\ly^5}h_{22;\,2}^2\right),
\umbruch\\
&\hspace*{-20ex}
\left(6\left(H^2-\A2\right)^2-12H^2\left(H^2-3\A2\right)\right)
\left(-F^{ij}H_{;\,i}H_{;\,j}\right)
=\umbruch\\
=&-\left(\lx^4+\lx^3\ly+\lx^2\ly^2+\lx\ly^3+\ly^4\right)\cdot\\
&\quad\cdot6\left(\lx^2+\ly^2\right)^2\frac{4}{\lx\ly}
\left(\frac1{\lx^5}h_{11;\,1}^2+\frac1{\ly^5}h_{22;\,2}^2\right)\umbruch\\
=&-\left(6\lx^8+6\lx^7\ly+18\lx^6\ly^2+18\lx^5\ly^3+24\lx^4\ly^4
+18\lx^3\ly^5\right.\\
&\qquad\left.+18\lx^2\ly^6+6\lx\ly^7+6\ly^8\right)\cdot\umbruch\\
&\quad\cdot\frac{4}{\lx\ly}
\left(\frac1{\lx^5}h_{11;\,1}^2+\frac1{\ly^5}h_{22;\,2}^2\right),\umbruch\\
-F^{ij}\Az i\Az j=&-\frac{4(\ly a_1+\lx)^2}{\lx^2\ly}h_{11;\,1}^2
-\frac{4(\lx a_2+\ly)^2}{\ly^2\lx}h_{22;\,2}^2\umbruch\\
=&-\left(\lx^3+\ly^3\right)^2\cdot\frac4{\lx\ly}
\left(\frac1{\lx^5}h_{11;\,1}^2+\frac1{\ly^5}h_{22;\,2}^2\right),\umbruch\\
2\left(H^2+2\A2\right)\cdot&\left(-F^{ij}\Az i\Az j\right)=\umbruch\\
=&-2\left(3\lx^2+2\lx\ly+3\ly^2\right)\left(\lx^3+\ly^3\right)^2\cdot\\
&\quad\cdot\frac4{\lx\ly}\left(\frac1{\lx^5}h_{11;\,1}^2
+\frac1{\ly^5}h_{22;\,2}^2\right)\umbruch\\
=&-\left(6\lx^8+4\lx^7\ly+6\lx^6\ly^2+12\lx^5\ly^3
+8\lx^4\ly^4+12\lx^3\ly^5\right.\\
&\qquad\left.+6\lx^2\ly^6+4\lx\ly^7+6\ly^8\right)
\frac4{\lx\ly}\left(\frac1{\lx^5}h_{11;\,1}^2
+\frac1{\ly^5}h_{22;\,2}^2\right),\umbruch\\
&\hspace*{-20ex}-F^{ij}\left(H_{;i}\Az j+H_{;j}\Az i\right)=\umbruch\\
=&-\frac{4(\ly a_1+\lx)(a_1+1)}{\lx^2\ly}h_{11;\,1}^2
-\frac{4(\lx a_2+\ly)(a_2+1)}{\ly^2\lx}h_{22;\,2}^2\umbruch\\
=&-\left(\lx^3+\ly^3\right)\left(\lx^2+\ly^2\right)
\cdot\frac4{\lx\ly}\left(\frac1{\lx^5}h_{11;\,1}^2
+\frac1{\ly^5}h_{22;\,2}^2\right),\umbruch\\
&\hspace*{-20ex}-12H\A2\left(-F^{ij}
\left(H_{;\,i}\Az j+H_{;\,j}\Az i\right)\right)=\umbruch\\
=&12(\lx+\ly)\left(\lx^2+\ly^2\right)^2
\left(\lx^3+\ly^3\right)\cdot\umbruch\\
&\quad\cdot\frac4{\lx\ly}
\left(\frac1{\lx^5}h_{11;\,1}^2
+\frac1{\ly^5}h_{22;\,2}^2\right)\umbruch\\
=&\left(12\lx^8+12\lx^7\ly+24\lx^6\ly^2
+36\lx^5\ly^3+24\lx^4\ly^4+36\lx^3\ly^5\right.\\
&\qquad\left.+24\lx^2\ly^6+12\lx\ly^7+12\ly^8\right)\cdot\\
&\quad\cdot\frac4{\lx\ly}\left(\frac1{\lx^5}h_{11;\,1}^2
+\frac1{\ly^5}h_{22;\,2}^2\right).
\end{align*}
Combining these expressions yields
\begin{align*}
&\hspace*{-20ex}\left(H^2-\A2\right)^4\left(\oprokl w\right)=\umbruch\\
=&-8(\lx+\ly)(\lx-\ly)^2\lx\ly
-\frac{32\ly^3}{\lx^2}h_{11;\,1}^2
-\frac{32\lx^3}{\ly^2}h_{22;\,2}^2,\umbruch\\
\oprokl w=&-\frac{(\lx+\ly)(\lx-\ly)^2}{2\lx^3\ly^3}
-\frac2{\lx^6\ly}h_{11;\,1}^2
-\frac2{\lx\ly^6}h_{22;\,2}^2\umbruch\\
\le&0.
\end{align*}
We finally apply the maximum principle and our theorem follows.
\end{proof}

\section{Convergence to Infinity}\label{inf sec}

It is known \cite{HuiskenPolden}, 
that \eqref{flow eqn} is a parabolic evolution equation
for strictly convex initial data and that it has a solution on a maximal
time interval $[0,\,T)$. Here, we want to show that 
some points on $M_t$ converge
to infinity for $t\uparrow T$, i.\,e.
$$\lim\limits_{t\uparrow T}\sup\limits_{M_t}|X|=\infty.$$

We show that the principal curvatures of $M_t$ stay uniformly 
bounded above.
\begin{lemma}\label{K upper bound}
For a smooth closed strictly convex surface $M$ in $\R^3$, 
flowing according to $\dot X=\frac1K\nu$, the maximum of the
principal curvatures is non-increasing.
\end{lemma}
\begin{proof}
Consider $M_{ij}=h_{ij}-\mu g_{ij}$ with $\mu>0$ so large
that $M_{ij}$ is negative semi-definite for some time $t_0$. We wish
to show that $M_{ij}$ is negative semi-definite for $t>t_0$. Combine
\eqref{h evol}, \eqref{Riem}, and \eqref{Ricci} to obtain
$$\dt h_{ij}-F^{kl}h_{ij;\,kl}=\frac1KH h_{ij}-\frac1Kh^k_ih_{kj}
-\frac1K\left(\tilde h^{kl}\tilde h^{rs}
+\tilde h^{kr}\tilde h^{ls}\right)h_{kl;\,i}h_{rs;\,j}.$$
In the evolution equation for $M_{ij}$, we drop the negative definite 
terms involving derivatives of the second fundamental form
$$\dt M_{ij}-F^{kl}M_{ij;\,kl}\le 
\frac1K\left(Hh_{ij}-h^k_ih_{kj}-2\mu h_{ij}\right).$$
Let $\xi$ be a zero eigenvalue of $M_{ij}$ with $\abs{\xi}=1$,
$M_{ij}\xi^j=h_{ij}\xi^j-\mu g_{ij}\xi^j=0$.
So we obtain in a point with $M_{ij}\le0$
$$\left(Hh_{ij}-h^k_ih_{kj}-2\mu h_{ij}\right)\xi^i\xi^j=
H\mu-3\mu^2\le2\mu\mu-3\mu^2\le0$$
and the maximum principle for tensors 
\cite{ChowLuMP,HamiltonThree} implies the claimed result.
\end{proof}

We obtain a pinching estimate
\begin{lemma}\label{pinching lem}
For a smooth closed strictly convex surface $M_t$ in $\R^3$, 
flowing according to $\dot X=\frac1K\nu$, there exists 
$c=c(M_0)$ such that $0<\frac1c\le\frac{\lambda_1}{\lambda_2}\le c$.
\end{lemma}
\begin{proof}
Choose $\mu>0$ such that $\lambda_1,\,\lambda_2\le\mu$
at $t=0$. Theorem \ref{K invers mon thm} and Lemma \ref{K upper bound}
imply that 
$$\frac1{\mu^2}\frac{\left(\frac\lx\ly-1\right)^2}{\frac\lx\ly}
=\frac1{\mu^2}\frac{(\lx-\ly)^2}{\lx\ly}\le
\frac{(\lx-\ly)^2}{\lx^2\ly^2}\le c.$$
We obtain the bound on $\frac{\lambda_1}{\lambda_2}$ claimed above. 
\end{proof}

It is only here that we use the monotone quantity of 
Theorem \ref{K invers mon thm}. For our purposes, this quantity
is better than scaling invariant. As it becomes apparent from 
the proof of Lemma \ref{1st C0 est}, however, the Aleksandrov
reflection principle can be used instead for the rest of the proof. 
This simplifies the
proof compared to \cite{AndrewsStones,OSA2}, where similar monotone
quantities are used. Later on, see Theorem \ref{impr K invers mon thm}
and Theorem \ref{impr conv thm}, 
we will use monotone quantities
to improve the convergence rate that follows from Aleksandrov
reflection. 

The next result shows that $K$ stays uniformly bounded below 
by a positive constant as
long as $M_t$ is enclosed by a ball of fixed positive radius. 
For similar results see \cite{TsoPoint} and 
\cite[Proposition 4.13]{ChowNonhomogGauss}.
\begin{lemma}\label{vel bound}
For a strictly convex solution of \eqref{flow eqn}, $K$ is
uniformly bounded below by a positive constant 
in terms of the radius $R$ of an enclosing 
sphere $B_{R}(x_0)$, the pinching ratio $\lx/\ly$, and
$\max_{M_0}\,\frac{K^{-1}}{2R-\langle X-x_0,\,\nu\rangle}$.
More precisely, we have everywhere
\begin{equation}\label{K bound}
K\ge\frac1{3R}\min\left\{\left(\max\limits_{M_0}
\frac{K^{-1}}{2R-\langle X-x_0,\,
\nu\rangle}\right)^{-1},\frac{\epsilon^2}{4R}\right\},
\end{equation}
where $\epsilon=\frac1{1+c}$ with $c$ as in Lemma \ref{pinching lem}. 
We obtain a positive lower bound on the principal curvatures. 
\end{lemma}
\begin{proof}
We may assume that $x_0=0$. Observe that
$3R\ge2R-\langle X,\,\nu\rangle\ge R$.
Thus $\frac{K^{-1}}{2R-\langle X,\,\nu\rangle}$ 
is finite for strictly convex surfaces. 
\par 
Standard computations \cite{HuiskenRoundSphere,OSMM,OSKnutCrelle,%
HuiskenPolden} 
yield the evolution equations
\begin{align*}
\dt X^\beta-F^{ij}X^\beta_{;\,ij}=&\frac3K\nu^\beta,\displaybreak[1]\\
\dt\nu^\beta-F^{ij}\nu^\beta_{;\,ij}=&\frac HK\nu^\beta,\displaybreak[1]\\
\dt\langle X,\,\nu\rangle-F^{ij}\langle X,\,\nu\rangle_{;\,ij}=&
-\frac1K+\frac HK\langle X,\,\nu\rangle,\displaybreak[1]\\
\opr{\frac1K}=&\frac H{K^2}.
\end{align*}
In a critical point of $\left.\frac{K^{-1}}{2R-\langle X,\,\nu\rangle}
\right|_{M_t}$, we obtain
$$\dt\log\frac{K^{-1}}{2R-\langle X,\,\nu\rangle}
-F^{ij}\left(\log\frac{K^{-1}}{2R-\langle X,\,\nu\rangle}\right)_{;\,ij}=
\frac1{2R-\langle X,\,\nu\rangle}\frac1K(2RH-1).$$
So we get in an increasing maximum
\begin{align*}
0\le&2RH-1,\\
\frac{\epsilon^2}{4R^2}\le&\epsilon^2H^2\le K
\end{align*}
and deduce there that
$$\frac{K^{-1}}{2R-\langle X,\,\nu\rangle}
\le\frac{4\epsilon^{-2}R^2}{R}.$$
Thus we obtain everywhere
$$\frac{K^{-1}}{3R}\le\frac{K^{-1}}{2R-\langle X,\,\nu\rangle}
\le\max\left\{\max\limits_{M_0}\frac{K^{-1}}{2R-\langle X,\,\nu\rangle},\,
\frac{4R}{\epsilon^2}\right\}$$
and \eqref{K bound} follows. 
\par
Finally, the positive lower bound on $K$ 
and our pinching estimate, Lemma \ref{pinching lem},
imply a positive lower bound on the principal curvatures. 
\end{proof}

Let us recall a form of the maximum principle for evolving hypersurfaces.
\begin{lemma}\label{encl lem}
Let $M_t$ and $\tilde M_t$ be two smooth closed strictly convex solutions
to \eqref{flow eqn} on some time interval $\left[0,\,T^*\right)$.
If $M_0$ encloses $\tilde M_0$, then $M_t$ encloses $\tilde M_t$
at any time $t\in\left[0,\,T^*\right)$,
for which both solutions exist.
\end{lemma}
\begin{proof}
This is a standard consequence of the maximum principle. 
\end{proof}

The next result describes the evolution of spheres. 
\begin{lemma}\label{sphere evol}
Spheres $\partial B_{r(t)}(x_0)$ solve \eqref{flow eqn} for
$t\in[0,\,T)$ with $r(t)=(T-t)^{-1}$ and $T=r^{-1}(0)$.
\end{lemma}
\begin{proof}
The evolution equation for the radius of a sphere is
$$\dot r(t)=r^2(t).$$
\end{proof}

\begin{lemma}
Let $M_t$ be a family of smooth closed strictly convex 
solutions to \eqref{flow eqn} on a maximal time interval
$[0,\,T)$. Then $T<\infty$.
\end{lemma}
\begin{proof}
Spheres solving \eqref{flow eqn} tend to infinity in 
finite time, see Lemma \ref{sphere evol}. So Lemma
\ref{encl lem} implies that $T<\infty$. 
\end{proof}

\begin{lemma}
Let $M_t$ be a family of smooth closed strictly convex
surfaces solving \eqref{flow eqn} on a maximal time interval
$[0,\,T)$. Then
$$\lim\limits_{t\uparrow T}\sup\limits_{M_t}\abs X=\infty.$$
\end{lemma}
\begin{proof}
Assume that $M_t\subset B_R(0)$, $0\le t<T$, for some $R>0$.
Then Lemmata \ref{K upper bound} and \ref{vel bound} imply
that the principal curvatures of $M_t$ stay uniformly bounded
above and below by positive constants. Thus \eqref{flow eqn}
can be rewritten as a uniformly parabolic equation and the 
estimates of Krylov, Safonov, Evans (see also
\cite{Andrews2dKrylov}), and Schauder imply uniform a priori
estimates up to $t=T$. This allows to continue the solution $M_t$
smoothly past $t=T$, a contradiction. Note finally, that
\eqref{flow eqn} is an expanding flow, so $\sup_{M_t}\abs X$
is monotone in $t$ and the limit exists. 
\end{proof}

\section{Convergence to a Sphere}\label{conv sphere sec}

\begin{lemma}\label{1st C0 est}
Under the assumptions of Theorem \ref{main thm}, 
we get
$$\lim\limits_{t\uparrow T}\inf\limits_{M_t}|X|=\infty,$$
more precisely, there exists $c=c(M_0)$ such that
$$\sup\limits_{M_t}|X|-c\le(T-t)^{-1}\le\inf\limits_{M_t}|X|+c$$
and
\begin{equation}\label{resc Hsd conv}
M_t\cdot (T-t)\subset B_{1+c(T-t)}(0)\setminus B_{1-c(T-t)}(0),
\end{equation}
so the rescaled surfaces $M_t\cdot(T-t)$ converge to the
unit sphere $\S^2$ in Hausdorff distance.
\end{lemma}
\begin{proof}
We may shift the origin such that $0$ lies inside $M_0$.
This does not affect the convergence rate claimed above.
\par
Define the support function $u:\S^2\times[0,\,T)\to\R_+$ for a 
convex surface by
$$u(z,\,t)=\langle X\left(\nu^{-1}(z),\,t\right),\,z\rangle.$$
It fulfills the evolution equation, see e.\,g.\ \cite{AndrewsHarnack},
$$\dot u=\frac{\det(u_{;\,ij}+u\sigma_{ij})}{\det(\sigma_{ij})},$$
where $u_{;\,ij}$ denotes covariant derivatives on $\S^2$
and $\sigma_{ij}$ is the standard metric on $\S^2$.
\par
We apply the Aleksandrov reflection principle of
Bennett Chow, Robert Gulliver \cite{ChowGulliverSphere},
and James McCoy \cite[Theorem 3.1]{JamesMcCoyArea} and
obtain a uniform bound on the oscillation (and the gradient)
of $u(\cdot,\,t)$ for all $t\in[0,\,T)$, that depends only
on the initial data. 
\par
As $\sup_{M_t}|X|\to\infty$ for $t\uparrow T$, we obtain that
$\inf_{M_t}|X|\to\infty$ for $t\uparrow T$, more precisely, we
have
$$\sup\limits_{M_t}|X|\le\inf\limits_{M_t}|X|+c.$$\par
It remains to show that 
$$\inf\limits_{M_t}|X|\le(T-t)^{-1}\le\sup\limits_{M_t}|X|.$$
Consider the surfaces $\partial M_{(T-t)^{-1}}(0)$, solving
\eqref{flow eqn}. For $t=T$, $M_t$ and 
$\partial B_{(T-t)^{-1}}(0)$ converge to infinity. We claim
that for all $t\in[0,\,T)$, $\partial B_{(T-t)^{-1}}(0)
\cap M_t\neq\emptyset$. Otherwise, for some $t_0$, $M_{t_0}$
encloses $\partial B_{(T-t)^{-1}}(0)$ or is contained in
$B_{(T-t)^{-1}}(0)$. Both cases are similar. We only
consider the first case. Choose $\epsilon>0$ such that
$M_{t_0}$ encloses also $\partial B_{(T-t_0-\epsilon)^{-1}}(0)$,
a slightly larger sphere. For $t\in[t_0,\,T-\epsilon)$,
$\partial B_{(T-t-\epsilon)^{-1}}(0)$ solves \eqref{flow eqn}.
As $\partial B_{(T-t-\epsilon)^{-1}}(0)$ converges to
infinity for $t\uparrow T-\epsilon$, Lemma \ref{encl lem}
implies that $M_t$ has to converge to infinity for
$t\uparrow T-\epsilon$, a contradiction.
\end{proof}

In terms of the support function $u$, this Lemma implies that
\begin{equation}\label{supp fct est}
u(x,\,t)-c\le(T-t)^{-1}\le u(x,\,t)+c.
\end{equation}
Note that this estimate is sharp for spheres
$\partial B_{(T-t)^{-1}}(Q)$ solving \eqref{flow eqn}, 
if $Q$ is different from the origin. The method of 
\cite{AndrewsStones}, see also \cite{OSA2}, where the
origin is replaced by some $q(t)$, can be adapted to the
present situation. In order to improve the estimate
\eqref{supp fct est}, however, we need a monotone
quantity similar to \eqref{mon groe} with a better
scaling behavior. We address this issue in Section 
\ref{impr conv sec}.

\section{Smooth Convergence to a Sphere}
\label{smooth conv sec}

\begin{lemma}\label{K upper bound two}
Under the assumptions of Theorem \ref{main thm}, there
exists $c=c(M_0)$, such that
$$K\le c\cdot(T-t)^2.$$
\end{lemma}
\begin{proof}
For $\mu\gg1$ to be fixed below, we consider
$$w:=|X|^2-\mu\frac1K.$$
We may assume that $\mu$ is so large that 
$w<0$ on $M_0$. Our aim is
to show that $w$ stays negative during the flow.
We use the evolution equations of the proof of Lemma \ref{vel bound}.
The evolution equation of $w$ is given by
$$\oprokl w=6\frac1K\langle X,\,\nu\rangle-2\frac H{K^2}
-\mu\frac H{K^2}.$$
Let $t_0\in[0,T)$ be minimal such that $\max_{M_t}w=0$.
Choose $x_0\in M_t$ such that $w(x_0,\,t_0)=0$. At this
point, we apply the parabolic maximum principle and 
obtain
\begin{align*}
0\le&6K\langle X,\,\nu\rangle-(2+\mu)H\umbruch\\
\le&6K|X|-(2+\mu)H\umbruch\\
=&6\sqrt K\sqrt\mu-(2+\mu)H&\text{as~}&w(x_0,\,t_0)=0,\umbruch\\
\le&(6\sqrt\mu-2-\mu)H&\text{as~}&K\le H^2.
\end{align*}
We fix $\mu$ sufficiently large and obtain that $w\le0$ during
the flow. In view of Lemma \ref{1st C0 est},
this implies the upper bound on the Gau{\ss} 
curvature claimed above.
\end{proof}

Combining this result with the Lemmata \ref{pinching lem},
\ref{vel bound}, and \ref{1st C0 est}, we obtain
\begin{equation}\label{K bounded}
\frac1c\cdot(T-t)^2\le K\le c\cdot(T-t)^2.
\end{equation}

We rescale our surfaces similarly as in \cite{AndrewsContractCalc}.
Consider the embeddings $\tilde X(\cdot,\,t)$,
$$\tilde X(z,\,t):=(T-t)\cdot X(z,\,t).$$
Define a new time function
$$\tau(t):=-\log\left(\frac{T-t}T\right).$$
We use a tilde to denote geometric quantities of the
rescaled surfaces. For $\tilde X$, we obtain the 
evolution equation
$$\frac d{d\tau}\tilde X=\frac1{\tilde K}\tilde\nu-\tilde X.$$
Our a priori estimates and the estimates of Krylov, Safonov,
Evans, and Schauder imply uniform bounds on all derivatives
of the support function $\tilde u$ of $\tilde X$. 
Applying interpolation inequalities as in
\cite[Lemma C.2]{OSKnutAIHP} to
\begin{align*}
\abs{\tilde u-1}\le&c\cdot(T-t),\umbruch\\
\abs{D^k\tilde u}\le&c_k,\umbruch
\intertext{we get}
\abs{D^k\tilde u}\le&c(k,\,\epsilon)\cdot(T-t)^{1-\epsilon}
\end{align*}
for any $\epsilon>0$. 

This finishes the proof of Theorem \ref{main thm}.

\section{Improved Convergence Rate}\label{impr conv sec}

\begin{theorem}\label{impr K invers mon thm}
For a family of smooth closed strictly convex surfaces $M_t$ 
in $\R^3$ flowing according to $\dot X=\frac1K\nu$, 
\begin{equation}
\max\limits_{M_t}\frac{\left(\lx^2+\ly^2\right)(\lambda_1-\lambda_2)^2}
{8(\lx+\ly)\lx^3\ly^3}
\equiv\max\limits_{M_t}w\equiv\max\limits_{M_t}
\frac{\A2\left(2\A2-H^2\right)}{H\left(H^2-\A2\right)^3}
\end{equation}
is non-increasing in time.
\end{theorem}
\begin{proof}
We use Section \ref{comp alg sec} and obtain in a critical point of $w$
\begin{align*}
\oprokl w=&\frac{-3\left(\lx^4+2\lx^3\ly-2\lx^2\ly^2+2\lx\ly^3+\ly^4\right)
(\lx-\ly)^2}{8(\lx+\ly)^2\lx^4\ly^4}\umbruch\\ 
&-\frac1{\left(3\lx^3+3\lx^2\ly-\lx\ly^2+3\ly^3\right)^2
(\lx+\ly)^2\lx^7\ly^2}\cdot\\
&\quad\cdot\left(18\lx^9-9\lx^8\ly+12\lx^7\ly^2
+72\lx^6\ly^3-12\lx^5\ly^4+70\lx^4\ly^5\right.\\
&\qquad\left.+60\lx^3\ly^6+18\lx\ly^8+27\ly^9\right)h_{11;\,1}^2\umbruch\\
&-\frac1{\left(3\ly^3+3\ly^2\lx-\ly\lx^2+3\lx^3\right)^2
(\ly+\lx)^2\ly^7\lx^2}\cdot\\
&\quad\cdot\left(18\ly^9-9\ly^8\lx+12\ly^7\lx^2
+72\ly^6\lx^3-12\ly^5\lx^4+70\ly^4\lx^5\right.\\
&\qquad\left.+60\ly^3\lx^6+18\ly\lx^8+27\lx^9\right)h_{22;\,2}^2\umbruch\\
\le&0.
\end{align*}
We finally apply the maximum principle. 
\end{proof}

This allows to improve our bound on $\abs{\lx-\ly}$.
\begin{lemma}\label{diff bound}
For a smooth closed strictly convex surface $M_t$ in $\R^3$,
flowing according to $\dot X=\frac1K\nu$, there exists a
constant $c=c(M_0)$ such that $\abs{\lx-\ly}\le c\cdot K^{5/4}
\le c\cdot (T-t)^{5/2}$.
\end{lemma}
\begin{proof}
This is a direct consequence of Theorem \ref{impr K invers mon thm},
Lemma \ref{pinching lem} and Lemma \ref{K upper bound two}.
\end{proof}

We now closely follow the corresponding parts of \cite{AndrewsStones}
and \cite{OSA2}.

\begin{proposition}\label{prop int H}
Define $q(t):=\frac1{4\pi}\int\limits_{M_t}KX$.
Then 
$$\left\lvert\langle X-q,\nu\rangle-\frac1{8\pi}
\int\limits_{M_t}H\right\rvert\le\frac1{4\pi}
\cdot\sup_{M_t}\abs{\lx-\ly}\cdot\mathcal H^{2}(M_t),$$
where $\mathcal H^2(M_t)$ denotes the area of $M_t$.
\end{proposition}
\begin{proof}
This is \cite[Proposition 4]{AndrewsStones}.
\end{proof}

We will call $q(t)$ the pseudocenter of $M_t$.

We define $r_+(t)$ to be the minimal radius of a sphere, 
centered at $q(t)$, that encloses $M_t$. Similarly, we
define $r_-(t)$ to be the maximal radius of a sphere,
centered at $q(t)$, that is enclosed by $M_t$.
Let $\rho_-(t)$ be the maximal radius of a sphere
(with arbitrary center) enclosed by $M_t$ and $\rho_+(t)$
be the minimal radius of spheres enclosing $M_t$.

\begin{lemma}\label{radii bound}
Under the assumptions of Theorem \ref{main thm}, 
for $T-t$ sufficiently small, $r_+$ and
$r_-$ are estimated as follows
$$(T-t)^{-1}\cdot\left(1-c\cdot(T-t)^{3/2}\right)\le r_-(t)
\le r_+(t)\le(T-t)^{-1}\cdot\left(1+c\cdot(T-t)^{3/2}\right).$$
\end{lemma}
\begin{proof}
Denote the bounded component of $\R^3\setminus M_t$
by $E_t$. The transformation formula for integrals
implies that
$$\frac1{4\pi}\int\limits_{M_t}KX=\frac1{4\pi}
\int\limits_{\S^2}X\left(\nu^{-1}(\cdot)\right).$$
So we see that $q(t)\in E_t$.
We have 
\begin{align*}
r_+=&\max_{M_t}\,\langle X-q(t),\,\nu\rangle,&
r_-=&\min_{M_t}\,\langle X-q(t),\,\nu\rangle,\\ 
\rho_+=&\min_{p\in\R^3}\max_{M_t}\,\langle X-p,\,\nu\rangle,
\qquad\text{and}&
\rho_-=&\max_{p\in E_t}\min_{M_t}\,\langle X-p,\,\nu\rangle.
\end{align*}

Recall the first variation formula for a vector 
field $Y$ along $M_t$
$$\int\limits_{M_t}H
\langle Y,\,\nu\rangle=\int\limits_{M_t}\divergenz_{M_t} Y$$
and get for $p\in E_t$ such that $\rho_+=\max_{M_t}
\langle X-p,\,\nu\rangle$
$$\int\limits_{M_t}H\ge\frac1{\rho_+}\int\limits_{M_t}
H\cdot\langle X-p,\,\nu\rangle=\frac1{\rho_+}\int\limits_{M_t}
\divergence_{M_t}X=\frac1{\rho_+}\int\limits_{M_t}2=\frac2{\rho_+}
\mathcal H^2(M_t).$$
We employ Proposition \ref{prop int H} and deduce that 
\begin{align*}
r_-\ge&\frac1{8\pi}\int\limits_{M_t}H\cdot\left\{1-
2\left(\int_{M_t}H\right)^{-1}\cdot
\sup\limits_{M_t}\abs{\lx-\ly}\cdot\mathcal H^2(M_t)\right\}\\
\ge&\frac1{8\pi}\int\limits_{M_t}H\cdot
\left\{1-\rho_+\cdot\sup\limits_{M_t}\abs{\lx-\ly}\right\}\\
\ge&\frac1{8\pi}\int\limits_{M_t}H\cdot
\left\{1-c\cdot(T-t)^{3/2}\right\},
\end{align*}
where we have used the Lemmata \ref{diff bound} and 
\ref{1st C0 est}.
So we obtain
\begin{equation}\label{r- bound}
r_-(t)\ge\frac1{8\pi}\int\limits_{M_t}H\cdot\left(1-
c\cdot(T-t)^{3/2}\right)
\end{equation}
and similarly
\begin{equation}\label{r+ bound}
r_+(t)\le\frac1{8\pi}\int\limits_{M_t}H\cdot\left(1+
c\cdot(T-t)^{3/2}\right).
\end{equation}
As in the proof of Lemma \ref{1st C0 est}, we get
$$r_-\le\rho_-\le(T-t)^{-1}\le\rho_+\le r_+$$
and 
$$(T-t)^{-1}\cdot\left(1-c\cdot(T-t)^{3/2}\right)
\le\frac1{8\pi}\int\limits_{M_t}H\le
(T-t)^{-1}\cdot\left(1+c\cdot(T-t)^{3/2}\right).$$
Using \eqref{r- bound} and \eqref{r+ bound} gives the
claimed estimates on $r_-$ and $r_+$.
\end{proof}

\begin{lemma}\label{Q est}
Under the assumptions of Theorem \ref{main thm}, 
$q(t)$ as defined in Proposition \ref{prop int H}
is a smooth function of $t$ in $[0,\,T)$ and converges
to some point $Q\in\R^3$ for $t\uparrow T$,
$$\abs{q(t)-Q}\le c\cdot (T-t)^{1/2}.$$
\end{lemma}
\begin{proof}
The definition of $q(t)$ involves only quantities that depend
smoothly on $t$, so it remains to prove convergence for
$t\uparrow T$.\par
For $0<t_1<t_2<T$, we want to estimate $\abs{q(t_1)-q(t_2)}$
from above. We may assume that $q(t_1)\neq q(t_2)$. Consider
the line passing through $q(t_1)$ and $q(t_2)$. It intersects
the surface $M_{t_2}$ in two points, denoted by $p_l(t_2)$ and
$p_r(t_2)$.
 
\setlength{\unitlength}{0.00075in}
\begingroup\makeatletter\ifx\SetFigFont\undefined%
\gdef\SetFigFont#1#2#3#4#5{%
  \reset@font\fontsize{#1}{#2pt}%
  \fontfamily{#3}\fontseries{#4}\fontshape{#5}%
  \selectfont}%
\fi\endgroup%
{\renewcommand{\dashlinestretch}{30}
\begin{picture}(6324,1135)(0,-10)
\drawline(12,748)(6312,748)
\drawline(3072,838)(3252,658)
\drawline(3252,838)(3072,658)
\drawline(2622,838)(2802,658)
\drawline(2802,838)(2622,658)
\drawline(552,1108)(551,1106)(549,1102)
	(546,1096)(541,1086)(536,1073)
	(529,1057)(521,1038)(513,1018)
	(504,996)(496,972)(489,947)
	(482,920)(476,891)(470,860)
	(466,825)(463,788)(462,748)
	(463,708)(466,671)(470,636)
	(476,605)(482,576)(489,549)
	(496,524)(504,500)(513,478)
	(521,458)(529,439)(536,423)
	(541,410)(546,400)(549,394)
	(551,390)(552,388)
\drawline(5772,1108)(5773,1106)(5775,1102)
	(5778,1096)(5783,1086)(5788,1073)
	(5795,1057)(5803,1038)(5811,1018)
	(5820,996)(5828,972)(5835,947)
	(5842,920)(5848,891)(5854,860)
	(5858,825)(5861,788)(5862,748)
	(5861,708)(5858,671)(5854,636)
	(5848,605)(5842,576)(5835,549)
	(5828,524)(5820,500)(5811,478)
	(5803,458)(5795,439)(5788,423)
	(5783,410)(5778,400)(5775,394)
	(5773,390)(5772,388)
\put(462,73){\makebox(0,0)[cb]{\smash{{\SetFigFont{12}{14.4}
{\rmdefault}{\mddefault}{\updefault}$p_l(t_2)$}}}}
\put(2712,73){\makebox(0,0)[cb]{\smash{{\SetFigFont{12}{14.4}
{\rmdefault}{\mddefault}{\updefault}$q(t_1)$}}}}
\put(3162,73){\makebox(0,0)[cb]{\smash{{\SetFigFont{12}{14.4}
{\rmdefault}{\mddefault}{\updefault}$q(t_2)$}}}}
\put(5862,73){\makebox(0,0)[cb]{\smash{{\SetFigFont{12}{14.4}
{\rmdefault}{\mddefault}{\updefault}$p_r(t_2)$}}}}
\end{picture}
}
\begin{center}Figure: Convergence of pseudocenters
\end{center}

We may assume that 
$$\langle p_r(t_2)-q(t_2),\,q(t_2)-q(t_1)\rangle>0.$$
This corresponds to $p_r(t_2)$ and $q(t_2)$ lying on the same
side of $q(t_1)$ as shown in the figure. We estimate
\begin{align*}
2c(T-t_1)^{1/2}\ge&r_+(t_1)-r_-(t_1) &&\text{by Lemma \ref{radii bound}}
\umbruch\\
=&\osc u_{q(t_1)}(\cdot,\,t_1)&&\umbruch\\
\ge&\osc u_{q(t_1)}(\cdot,\,t_2)&&
\text{by \cite[Theorem 3.1]{JamesMcCoyArea}}\umbruch\\
\ge&\abs{p_r(t_2)-q(t_1)}-\abs{p_l(t_2)-q(t_1)}&&\umbruch\\
=&\abs{p_r(t_2)-q(t_2)}+\abs{q(t_2)-q(t_1)}&&\\
&-\left(\abs{p_l(t_2)-q(t_2)}-\abs{q(t_2)-q(t_1)}\right)&&\umbruch\\
\ge&2\abs{q(t_2)-q(t_1)}+r_-(t_2)-r_+(t_2)&&\umbruch\\
\ge&2\abs{q(t_2)-q(t_1)}-2c(T-t_2)^{1/2},&&
\end{align*}
where $u_{q(t_1)}(\cdot,\,t_2)$ is the support function of
$M_{t_2}-q(t_1)$. Thus we can apply Cauchy's convergence criterion.
Finally, we let $t_2\uparrow T$ and the claimed bound follows.
\end{proof}

This allows to improve our convergence result.
\begin{theorem}\label{impr conv thm}
Under the assumptions of Theorem \ref{main thm}, there exists
$Q\in\R^3$ such that the Hausdorff distance of $M_t$ to a
family of 
expanding spheres around $Q$ is bounded as follows
$$d_{\mathcal H}(M_t,\,\partial B_{(T-t)^{-1}}(Q))\le
c\cdot(T-t)^{1/2}.$$
\end{theorem}
\begin{proof}
Combine the Lemmata \ref{Q est} and \ref{radii bound}.
\end{proof}

\begin{remark}
Theorem \ref{impr conv thm} implies also better estimates for
rescaled surfaces,
$$d_{\mathcal H}\left((M_t-Q)\cdot(T-t),\,\S^2\right)\le
c\cdot(T-t)^{3/2}.$$
As above, this implies that all derivatives of the support
function decay,
$$\Vert u_Q(\cdot,\,t)-1\Vert_{C^k}\le c(k,\,\epsilon)
\cdot(T-t)^{3/2-\epsilon}$$
for any $\epsilon>0$, where $u_Q(\cdot,\,t)$ is the support
function of $M_t-Q$. This implies for the principal 
curvatures $\tilde\lambda_i$, $i=1,\,2$, of $(M_t-Q)\cdot
(T-t)$
$$1-c(\epsilon)\cdot(T-t)^{3/2-\epsilon}\le
\tilde\lambda_i\le1+c(\epsilon)\cdot(T-t)^{3/2-\epsilon}.$$
\par
Without the additional $Q$, we get for expanding spheres
$\partial B_{(T-t)^{-1}}(P)$, 
$$d_{\mathcal H}\left(\partial B_{(T-t)^{-1}}(P)\cdot(T-t),\,\S^2\right)
=\abs P\cdot(T-t)^{-1},$$
so the estimate in Lemma \ref{resc Hsd conv} is sharp for
$P\neq0$.
\par
In the proof of Theorem \ref{main thm}, we have used 
\eqref{mon groe} only to prove that surfaces stay uniformly
pinched, i.\,e.\ that $\lambda_1/\lambda_2$ is uniformly bounded.
If we use it to bound $\abs{\lx-\ly}$, and then $r_+$ and $r_-$
as above, we don't get more than in \ref{resc Hsd conv}, where
we used the oscillation estimates of \cite{JamesMcCoyArea}. 
Our computer program, however, did not yield a scaling invariant
quantity that implies uniform pinching. This is similar to
\cite{Andrews2dnonconcave}. It might be possible to find a 
monotone quantity that allows to further improve this 
convergence rate. 
\end{remark}


\section{Finding Monotone Quantities}\label{algor sec}

\subsection{The Algorithm}
We use a sieve algorithm and start with symmetric rational functions 
of the principal curvatures
as candidates for test functions, e.\,g.\ 
$$w=\frac{p_1(\lx,\,\ly)}{p_2(\lx,\,\ly)}
=\frac{(\lx-\ly)^2}{\lx^2\ly^2}.$$
Here, $p_1\neq0$ and $p_2\neq0$ are homogeneous polynomials. 

In the end, we want to find functions $w$ such that 
$W:=\sup_{M_t}w$ is monotone and ensures convergence to round
spheres. 

We check, whether these test functions $w$ fulfill the
following conditions. 
\begin{enumerate}
\item\label{step one} 
 \begin{enumerate}
 \item $p_1(\lx,\,\ly),\,p_2(\lx,\,\ly)\ge0$ for $0<\lx,\,\ly$,
 \item $p_1(\lx,\,\ly)=0$ for $\lx=\ly>0$.
 \end{enumerate}
\item\label{step two} $\deg p_1<\deg p_2$.
\item\label{step three} $\fracp{w(1,\ly)}{\ly}<0$ for $0<\ly<1$ and 
  $\fracp{w(1,\ly)}{\ly}>0$ for $\ly>1$.
\item\label{step four} $\dt w-F^{ij}w_{;\,ij}\le0$ 
 \begin{enumerate}
 \item\label{step four a} for terms without derivatives of $(h_{ij})$,
 \item\label{step four b} for terms involving derivatives of $(h_{ij})$,
   if $w_{;\,i}=0$ for $i=1,\,2$.
 \end{enumerate}
\end{enumerate}

\subsection{Motivation and Randomized Tests}
For all flow equations 
considered, spheres contract to points and stay spherical. So
we can only find monotone quantities, if $\deg p_1\ge\deg p_2$
or $p_1(\lambda,\,\lambda)=0$. 

If $\deg p_1\ge\deg p_2$,
we obtain that $W$ is non-increasing on any self-similarly
expanding surface. So this does not imply convergence to a
sphere. 

Condition \eqref{step three} ensures that the quantity decreases, if 
the eigenvalues approach each other. 

In step \eqref{step four a} and \eqref{step four b}, we check that 
we can apply the maximum principle. Here, we have to use various
differentiation rules.

In steps \eqref{step one}, \eqref{step two}, and \eqref{step three}, 
inequalities are tested by evaluating both 
sides at random numbers. After enough testing, all candidates
for which the above inequalities, evaluated at random numbers,
were not violated, could be used to prove convergence to 
round spheres at infinity.

Alternatively, for surfaces, we can avoid using random numbers, 
compute evolution equations algebraically, and use Sturm's 
algorithm to test for non-negativity. 

We expect that similar algorithms will be used to find monotone
test functions for other (geometric) problems. 

\section{Computing Evolution Equations}\label{comp alg sec}

It is straightforward to use a computer algebra system 
to obtain evolution equations of test quantities $w$,
evaluated at a critical point of $w$. More precisely, let
$M_t$ be a family of surfaces in $\R^3$, moving with 
normal velocity $F=F(\lx,\,\ly)$, where $F>0$ for 
contracting surfaces. Assume that the test quantity $w$
is a function of $H$ and $\A2$. Then $w$ fulfills the
evolution equation
$$\opr w=C_w(\lx,\,\ly)+G_w(\lx,\,\ly)h_{11;\,1}^2
+G_w(\ly,\,\lx)h_{22;\,2}^2$$
in a critical point of $w$. It remains to 
compute $C_w$ (``constant terms'') and $G_w$ (``gradient terms''). 
The following calculations are all
similar as before and use \eqref{g evol}, \eqref{h evol},
\eqref{Riem}, \eqref{Ricci}, and, see 
\cite{AndrewsContractCalc,CGJDG1996},
$$F^{ij,\,kl}\eta_{ij}\eta_{kl}=\sum\fracp{^2F}{\li\partial\lj}
\eta_{ii}\eta_{jj}+\sum\limits_{i\neq j}
\frac{\fracp F\li-\fracp F\lj}{\li-\lj}(\eta_{ij})^2,$$
for symmetric matrices $(\eta_{ij})$ and $\lx\neq\ly$
or $\lx=\ly$ and the last term is interpreted as a limit.
For $w=H$, we obtain
\begin{align*}
C_H=&\left(\fracp F\lx \lx^2+\fracp F\ly \ly^2\right)H
+\left(F-\fracp F\lx\lx-\fracp F\ly\ly\right)\A2,\umbruch\\
G_H=&\fracp{^2F}{\lx\partial\lx}+2\fracp{^2F}{\lx\partial\ly}a_1
+\fracp{^2F}{\ly\partial\ly}a_1^2
+2\frac{\fracp F\lx-\fracp F\ly}{\lx-\ly}a_1^2,
\end{align*}
where
$$a_1=-\fracp{w\left(\lx+\ly,\,\lx^2+\ly^2\right)}{\lx}
\left(\fracp{w\left(\lx+\ly,\,\lx^2+\ly^2\right)}{\ly}\right)^{-1}$$
is such that $h_{22;\,1}=a_1h_{11;\,1}$ in a critical point
of $w$.
Similarly, we get for $w=\A2$
\begin{align*}
C_{\A2}=&2\left(\fracp F\lx\lx^2+\fracp F\ly\ly^2\right)\A2
+2\left(F-\fracp F\lx\lx-\fracp F\ly\ly\right)\A3,\umbruch\\
G_{\A2}=&-2\left(\fracp F\lx\left(1+a_1^2\right)+2\fracp F\ly a_1^2\right)
\umbruch\\
&+2\left(\fracp{^2F}{\lx\partial\lx}+2\fracp{^2F}{\lx\partial\ly}a_1
+\fracp{^2F}{\ly\partial\ly}a_1^2\right)\lx
+4\frac{\fracp F\lx-\fracp F\ly}{\lx-\ly}a_1^2\ly.
\end{align*}
We also need some mixed terms
\begin{align*}
-F^{ij}H_{;i}H_{;j}=&-\fracp F\lx(1+a_1)^2h_{11;\,1}^2
-\fracp F\ly(1+a_2)^2h_{22;\,2}^2,\umbruch\\
-F^{ij}\Az i\Az j=&-4\fracp F\lx(\lx+a_1\ly)^2h_{11;\,1}^2
-4\fracp F\ly(\ly+a_2\lx)^2h_{22;\,2}^2,\umbruch\\
-2F^{ij}H_{;i}\Az j=&-F^{ij}
\left(H_{;i}\Az j+H_{;j}\Az i\right)\umbruch\\
=&-4\fracp F\lx(1+a_1)(\lx+a_1\ly)h_{11;\,1}^2\\
&-4\fracp F\ly(1+a_2)(\ly+a_2\lx)h_{22;\,2}^2.
\end{align*}
Combining these expressions yields
\begin{align*}
C_w=&\fracp wHC_H+\fracp w{\A2}C_{\A2},\umbruch\\
G_w=&\fracp wHG_H+\fracp w{\A2}G_{\A2}
-\fracp{^2w}{H\partial H}\fracp F\lx(1+a_1)^2\\
&-4\fracp{^2w}{\A2\partial\A2}\fracp F\lx(\lx+a_1\ly)^2
-4\fracp{^2w}{H\partial\A2}\fracp F\lx(1+a_1)(\lx+a_1\ly).
\end{align*}
This formulae allow to easily compute evolution equations
in critical points.

\section{Other Normal Velocities}\label{o nor vel}

\subsection{Homogeneity less than minus one}

In this section, we prove 
\begin{theorem}\label{sammel thm}
For any smooth closed strictly convex surface $M$ in $\R^3$, 
there exists a smooth family of surfaces $M_t$, $t\in[0,T)$,
solving one of the following flow equations
\begin{itemize}
\item ${\displaystyle\dt X=\frac{H^2}{K^2}\nu,}$\\[0.5ex]\rule{1mm}{0mm}
\item ${\displaystyle\dt X=\frac{\A2}{K^2}\nu,}$\\[0.5ex]\rule{1mm}{0mm}
\item ${\displaystyle\dt X=\frac{H^3}{K^3}\nu,}$\\[0.5ex]\rule{1mm}{0mm}
\end{itemize}
with $M_0=M$.
For $t\uparrow T$, $M_t$ converges to infinity. 
The rescaled surfaces $M_t\cdot r^{-1}(t)$ converge smoothly to the
unit sphere $\S^2$, where $r(t)$ is the radius of an expanding
sphere that converges to infinity for $t\uparrow T$, more
precisely, $r(t)$ is as follows
\begin{itemize}
\item ${\displaystyle\dt X=\frac{H^2}{K^2}\nu:
\quad r(t)=(4(T-t))^{-1},}$\\[0.5ex]\rule{1mm}{0mm}
\item ${\displaystyle\dt X=\frac{\A2}{K^2}\nu:
\quad r(t)=(2(T-t))^{-1},}$\\[0.5ex]\rule{1mm}{0mm}
\item ${\displaystyle\dt X=\frac{H^3}{K^3}\nu:
\quad r(t)=(16(T-t))^{-1/2}.}$
\end{itemize}
\end{theorem}

We have also studied convex surfaces contracting according
to $$\dt X=-F\nu,$$
where $F$ is positive homogeneous of some degree larger than 
or equal to two \cite{OSA2}. 
There, we got the impression, that appropriate monotone quantities
that assure convergence to a sphere after rescaling are available
for almost every normal velocity considered.  
In contrast to this, such monotone quantities seem to be rare 
objects for expanding surfaces with normal velocity of homogeneity
less than or equal to minus two. 
\par
In both cases, we restricted our attention to normal velocities
and possible candidates for monotone quantities in
$\Z(\lx,\,\ly)$ with small coefficients which are 
symmetric in $\lx$ and $\ly$. 
\par
In the case of expanding surfaces, however, we also find
monotone quantities for surfaces expanding with a normal 
velocity that is positive homogeneous of degree minus one in 
the principal curvatures. 

The proof of Theorem \ref{sammel thm} is similar to the proof
of Theorem \ref{main thm}. Therefore, we will present in the
following only those results that are not almost identical 
to their counterparts in the proof of Theorem \ref{main thm}.

\begin{theorem}
For a family $M_t$ of smooth closed strictly convex surfaces in 
$\R^3$, flowing according to $\dt X=\frac{H^2}{K^2}\nu$, 
$$\max_{M_t}\frac{(\lx-\ly)^2}{2(\lx+\ly)\lx\ly}
\equiv\max_{M_t}w$$
is non-increasing in time.
\end{theorem}
\begin{proof}
According to Section \ref{comp alg sec}, we obtain in a
critical point of $w$
\begin{align*}
\oprokl w=&-\frac{5(\lx-\ly)^2}{\lx^2\ly^2}\\
&-\frac{128}{(\lx+3\ly)^2\lx^4}\cdot h_{11;\,1}^2
+(\ldots)\cdot h_{22;\,2}^2
\end{align*}
and apply the maximum principle. 
\end{proof}

The factor $2$ in the denominator is useful to rewrite $w$
in terms of the algebraic basis consisting of $H$ and $\A2$
$$\frac{(\lx-\ly)^2}{2(\lx+\ly)\lx\ly}=
\frac{2\A2-H^2}{H\cdot\left(H^2-\A2\right)}.$$

\begin{theorem}
For a family $M_t$ of smooth closed strictly convex surfaces in 
$\R^3$, flowing according to $\dt X=\frac{\A2}{K^2}\nu$, 
$$\max_{M_t}\frac{(\lx-\ly)^2}{2(\lx+\ly)\lx\ly}
\equiv\max_{M_t}w$$
is non-increasing in time.
\end{theorem}
\begin{proof}
According to Section \ref{comp alg sec}, we obtain in a
critical point of $w$
\begin{align*}
\oprokl w=&-\frac{2\left(2\lx^2+\lx\ly+2\ly^2\right)(\lx-\ly)^2}
{(\lx+\ly)^2\lx^2\ly^2}\\ 
&-\frac{4\left(21\lx^4+24\lx^3\ly+18\lx^2\ly^2+\ly^4\right)}
{(\lx+3\ly)^2(\lx+\ly)^2\lx^6}
\cdot h_{11;\,1}^2
+(\ldots)\cdot h_{22;\,2}^2
\end{align*}
and apply the maximum principle. 
\end{proof}

\begin{theorem}
For a family $M_t$ of smooth closed strictly convex surfaces in 
$\R^3$, flowing according to $\dt X=\frac{H^3}{K^3}\nu$, 
$$\max_{M_t}\frac{(\lx+\ly)^6(\lx-\ly)^2}{16\left(\lx^2+\ly^2\right)
\left(\lx^2+\lx\ly+\ly^2\right)\lx^3\ly^3}
\equiv\max_{M_t}w$$
is non-increasing in time.
\end{theorem}
\begin{proof}
According to Section \ref{comp alg sec}, we obtain in a
critical point of $w$
\begin{align*}
\oprokl w=&-\frac{(\lx+\ly)^8(\lx-\ly)^2}{4\left(\lx^2+\ly^2\right)^2
\left(\lx^2+\lx\ly+\ly^2\right)^2\lx^6\ly^6}\cdot\\
&\quad\cdot(\lx^6-\lx^5\ly+8\lx^4\ly^2+2\lx^3\ly^3+8\lx^2\ly^4
-\lx\ly^5+\ly^6)\\ 
&-\frac{3(\lx+\ly)^8}{8\left(3\lx^6+11\lx^4\ly^2+2\lx^3\ly^3
+9\lx^2\ly^4-2\lx\ly^5+\ly^6\right)^2}
\cdot\\
&\quad\cdot\frac1{\left(\lx^2+\lx\ly+\ly^2\right)^2
\left(\lx^2+\ly^2\right)^2\lx^8\ly^6}\cdot\\
&\quad\cdot\left(2\lx^{18}-4\lx^{17}\ly+57\lx^{16}\ly^2
-108\lx^{15}\ly^3+508\lx^{14}\ly^4\right.\\
&\qquad\left.-428\lx^{13}\ly^5+2152\lx^{12}\ly^6-156\lx^{11}\ly^7
+4784\lx^{10}\ly^8\right.\\
&\qquad\left.+172\lx^9\ly^9+4942\lx^8\ly^{10}-612\lx^7\ly^{11}
+2676\lx^6\ly^{12}\right.\\
&\qquad\left.-772\lx^5\ly^{13}+872\lx^4\ly^{14}-340\lx^3\ly^{15}
+126\lx^2\ly^{16}\right.\\
&\qquad\left.-56\lx\ly^{17}+9\ly^{18}\right)\cdot h_{11;\,1}^2\\
&+(\ldots)\cdot h_{22;\,2}^2.
\end{align*}
Here and in the following, we write $(\ldots)$ to denote a 
term that equals the factor in front of $h_{11;\,1}^2$
with $\lx$ and $\ly$ interchanged. 

We finally apply the maximum principle. 
\end{proof}

Similarly as above, we obtain the following evolution equations
for a family $M_t$ of surfaces flowing according to $\dt X=-F\nu$
\begin{align*}
\dt X^\alpha-F^{ij}X^\alpha_{;\,ij}
=&\left(F^{ij}h_{ij}-F\right)\nu^\alpha,\umbruch\\
\opr{\abs{X}^2}=&2\left(F^{ij}h_{ij}-F\right)\langle X,\,\nu\rangle
-2F^{ij}g_{ij},\umbruch\\
\opr{\nu^\alpha}=&F^{ij}h^k_ih_{kj}\cdot\nu^\alpha,\umbruch\\
\oprokl{\langle X,\,\nu\rangle}=&-F^{ij}h_{ij}-F+F^{ij}h^k_ih_{kj}
\langle X,\,\nu\rangle,\umbruch\\
\oprokl F=&FF^{ij}h^k_ih_{kj},\umbruch\\
\dt g_{ij}=&-2Fh_{ij},\umbruch\\
\dt h_{ij}-F^{kl}h_{ij;\,kl}=&F^{kl}h^a_kh_{al}\cdot h_{ij}
-F^{kl}h_{kl}\cdot h^a_ih_{aj}\\
&-Fh^k_ih_{kj}+F^{kl,\,rs}h_{kl;\,i}h_{rs;\,j}.
\end{align*}

\begin{lemma}\label{upper princ curv bnd}
Under the assumptions of Theorem \ref{sammel thm},
$\dt X=-F\nu$, the maximum principal curvature of $M_t$
is non-increasing in time.
\end{lemma}
\begin{proof}
Consider $M_{ij}:=h_{ij}-\mu g_{ij}$, where $\mu$ is chosen
such that $M_{ij}\le0$ initially. We compute
\begin{align*}
\dt M_{ij}-F^{kl}M_{ij;\,kl}=&
F^{kl}h^a_kh_{al}\cdot h_{ij}-F^{kl}h_{kl}\cdot h^a_ih_{aj}
-Fh^k_ih_{kj}+2\mu Fh_{ij}\\
&+F^{kl,\,rs}h_{kl;\,i}\xi^ih_{rs;\,j}\xi^j.
\end{align*}
Assume that $\xi$ is a zero eigenvalue of $(M_{ij})$,
$h_{ij}\xi^j=\mu g_{ij}\xi^j$, with $g_{ij}\xi^i\xi^j=1$. 
We may assume, that in our coordinate system, we have
$\xi=(1,\,0)$ and $(h_{ij})=\left(\begin{smallmatrix}
\mu&0\\ 0&\lambda\end{smallmatrix}\right)$ with $0<\lambda\le\mu$.
Normal velocities $F$ as in Theorem \ref{sammel thm} are
concave, so
$$F^{kl,\,rs}h_{kl;\,i}\xi^ih_{rs;\,j}\xi^j\le0.$$
We estimate
\begin{itemize}
\item $\left(\dt M_{ij}-F^{kl}M_{ij;\,kl}\right)\xi^i\xi^j\le
-\frac{(\lambda+\mu)(3\mu-\lambda)}
{\lambda^2}\le0$ for $F=-\frac{H^2}{K^2}$,
\item $\left(\dt M_{ij}-F^{kl}M_{ij;\,kl}\right)\xi^i\xi^j\le
-\frac{-2\lambda\mu+\lambda^2+3\mu^2}{\lambda^2}\le0$
for $F=-\frac{\A2}{K^2}$,
\item $\left(\dt M_{ij}-F^{kl}M_{ij;\,kl}\right)\xi^i\xi^j\le
-\frac{2(\lambda+\mu)^2(2\mu-\lambda)}{\mu\lambda^3}\le0$
for $F=-\frac{H^3}{K^3}$.
\end{itemize}
Now the claim follows directly from the maximum principle
for tensors.
\end{proof}

\begin{lemma}
Under the assumptions of Theorem \ref{sammel thm}, we obtain
a lower bound on the principal curvatures similarly to
Lemma \ref{vel bound}.
\end{lemma}
\begin{proof}
We proceed as in the proof of Lemma \ref{vel bound}.
In a critical point of $\frac{-F}{2R-\langle X,\,\nu\rangle}$, 
we compute the following evolution equation for a family
of surfaces flowing according to $\dt X=-F\nu$
\begin{align*}
\dt\log\frac{-F}{2R-\langle X,\,\nu\rangle}
-F^{ij}\left(\log\frac{-F}
{2R-\langle X,\,\nu\rangle}\right)_{;\,ij}=&
\frac1{2R-\langle X,\,\nu\rangle}\cdot\\
&\cdot\left(2R\cdot F^{ij}h^k_ih_{kj}-F^{ij}h_{ij}-F\right).
\end{align*}
We compute this explicitly for $F$ as in Theorem \ref{sammel thm}
and obtain in an increasing maximum of 
$\frac{-F}{2R-\langle X,\,\nu\rangle}$
\begin{itemize}
\item $8R-\frac HK\ge0$ for $F=-\frac{H^2}{K^2}$,
\item $4HR-\frac{\A2}K\ge0$ for $F=-\frac{\A2}{K}$,
\item $6R-\frac HK\ge0$ for $F=-\frac{H^3}{K^3}$.
\end{itemize}
Our monotone quantities and Lemma \ref{upper princ curv bnd}
imply that our surfaces are pinched, i.\,e.\ that 
$\frac\lx\ly+\frac\ly\lx$ is uniformly bounded above. So
we obtain there that $\lambda_i\ge\frac1{cR}$.
\end{proof}

\begin{lemma}
Consider one of the flow equations of Theorem \ref{sammel thm}, 
$\dt X=-F\nu$, and a solution as in this theorem.
Let $\gamma$ be such that $F$ is positive
homogeneous of degree $-\gamma$. Then there exists a
constant $\alpha\gg1$ such that $\abs F\ge\frac1\alpha
|X|^\gamma$ during the flow. More precisely, there
exists $\alpha\gg1$ such that $|X|^\gamma+\alpha F$ remains
non-positive during the flow, if this quantity is negative 
initially.  
\end{lemma}
\begin{proof}
For $F$ positive homogeneous of degree minus two, 
we get $F^{ij}h_{ij}=-2F$ and obtain at a point,
where $|X|^2+\alpha F=0$,
\begin{align*}
\dt\left(|X|^2+\alpha F\right)
-F^{ij}\left(|X|^2+\alpha F\right)_{;\,ij}=&
-6F\langle X,\,\nu\rangle-2F^{ij}g_{ij}
+\alpha FF^{ij}h^k_ih_{kj}\umbruch\\
\le&-6F\sqrt\alpha\sqrt{-F}
+\alpha FF^{ij}h^k_ih_{kj}.
\end{align*}
It is straightforward to check that for $\alpha\gg1$ sufficiently
large and $F$ as in the lemma, 
the right-hand side is non-positive. 
\par
If $F=-\frac{H^3}{K^3}$, we obtain similarly as above
\begin{align*}
\dt\left(|X|^3+\alpha F\right)
-F^{ij}\left(|X|^3+\alpha F\right)_{;\,ij}\le&
12|X|^2\frac{H^3}{K^3}-6\alpha\frac{H^5}{K^5}\le 0.
\end{align*}
\par
In both cases, the lemma follows from the maximum principle.
\end{proof}

\subsection{Homogeneity minus one}

\begin{table}
\def\platz{\raisebox{0em}[2.2em][1.5em]{\rule{0em}{2em}}}
$$\begin{array}{|c||c|}\hline
-\fracd1H & \platz\fracd{(\lx-\ly)^2}{(\lx+\ly)\lx\ly}\\\hline
-\fracd1H & \platz\fracd{\left(\lx^2+\ly^2\right)(\lx-\ly)^2}
{(\lx+\ly)\lx^3\ly^3}\\\hline
-\fracd HK & \platz\fracd{(\lx-\ly)^2}{\lx^2\ly^2}\\\hline
\end{array}$$
\caption{More monotone quantities}\label{tab2}
\end{table}

In Table \ref{tab2}, we have collected some normal velocities $F$,
positive homogeneous of degree minus one, and test functions $w$,
such that $\max_{M_t}w$ is non-increasing during the flow of
a closed strictly convex surface flowing according to
$$\dt X=-F\nu.$$

\begin{theorem}
For a family $M_t$ of smooth closed strictly convex surfaces in 
$\R^3$, flowing according to $\dt X=\frac1H\nu$, 
$$\max_{M_t}\frac{(\lx-\ly)^2}{2(\lx+\ly)\lx\ly}
\equiv\max_{M_t}w$$
and 
$$\max_{M_t}
\frac{\left(\lx^2+\ly^2\right)(\lx-\ly)^2}{8(\lx+\ly)\lx^3\ly^3}
\equiv\max_{M_t}\hat w$$
are non-increasing in time.
\end{theorem}
\begin{proof}
According to Section \ref{comp alg sec}, we obtain in a
critical point of $w$
\begin{align*}
\oprokl w=&-\frac{\left(\lx^2+4\lx\ly+\ly^2\right)(\lx-\ly)^2}
{2(\lx+\ly)^3\lx\ly}\\
&-\frac{2\left(5\lx^2+2\lx\ly+\ly^2\right)\ly}
{(\lx+3\ly)^2(\lx+\ly)\lx^5}\cdot h_{11;\,1}^2\\
&+(\ldots)\cdot h_{22;\,2}^2
\end{align*}
and apply the maximum principle. 
\par
Similarly, we obtain in a critical point of $\hat w$
\begin{align*}
\oprokl{\hat w}=&-\frac{\left(3\lx^4-2\lx^2\ly^2+3\ly^4\right)(\lx-\ly)^2}
{8(\lx+\ly)^3\lx^3\ly^3}\umbruch\\
&-\frac1{2\left(3\lx^3+3\lx^2\ly-\lx\ly^2+3\ly^3\right)^2
(\lx+\ly)^3\lx^7\ly}\cdot\\
&\quad\cdot\left(9\lx^{10}-9\lx^8\ly^2+96\lx^7\ly^3
-38\lx^6\ly^4+96\lx^5\ly^5+30\lx^4\ly^6\right.\\
&\qquad\left.+45\lx^2\ly^8+27\ly^{10}\right)\cdot h_{11;\,1}^2\\
&+(\ldots)\cdot h_{22;\,2}^2
\end{align*}
and apply the maximum principle once again. 
\end{proof}

\begin{theorem}
For a family $M_t$ of smooth closed strictly convex surfaces in 
$\R^3$, flowing according to $\dt X=\frac HK\nu$, 
$$\max_{M_t}\frac{(\lx-\ly)^2}{4\lx^2\ly^2}
\equiv\max_{M_t}w$$
is non-increasing in time.
\end{theorem}
\begin{proof}
According to Section \ref{comp alg sec}, we obtain in a
critical point of $w$
\begin{align*}
\oprokl w=&-\frac{(\lx-\ly)^2}{\lx^2\ly^2}
-\frac2{\lx^6}\cdot h_{11;\,1}^2+(\ldots)\cdot h_{22;\,2}.
\end{align*}
and apply the maximum principle. 
\end{proof}

For these flow equations, convergence to infinity and
convergence to a sphere after rescaling have been proved
before for hypersurfaces 
\cite{CGFlowSpheres,UrbasExpandMZ,UrbasExpandJDG,KnutHarmonicMCF}. 
Monotone quantities as mentioned above might at most be useful
to improve the convergence rate. Inverse mean curvature
was used to prove the Penrose inequality in general 
relativity \cite{HuiskenIlmanenPenrose}. Our techniques
might also apply to surfaces expanding in the asymptotically 
flat manifolds considered there.

\section{Convergence Rate}\label{conv rate}

In order to find out what the optimal convergence rate might
be, we proceed as in \cite{OSA2}, use the same notation, 
and compute the linearized equation corresponding to
$$\dt X=\frac1K\nu-X\qquad\text{or}\qquad
\fracp ut=\frac1K w-u$$
as
$$\fracp vt=\Delta v+v.$$
As in the contracting case, we only need to consider 
eigenvalues $-l(l+1)$ of the laplacian on the sphere 
for $l\in\N_+$. For $l=1$, we cannot expect convergence 
rates better than
$$r_+\le(T-t)^{-1}\cdot(1+c\cdot(T-t)),$$
if we fix $q(t)$ arbitrarily. This estimate is sharp,
if we do not adjust $q(t)$. The corresponding eigenfunctions
induce translations of the surface. Considering $l=2$, we 
expect that we cannot obtain convergence rates better than
$$r_+(t)\le(T-t)^{-1}\cdot\left(1+c\cdot(T-t)^5\right).$$

\bibliographystyle{/home/napier/schnuere/os/tools/amsplain}

\begin{thebibliography}{10}

\bibitem{Andrews2dKrylov}
Ben{}{{}}{} Andrews, \emph{{Fully nonlinear parabolic equations in two space
  variables}}, \relax\\{\tt arXiv:math.AP/0402235}.

\bibitem{Andrews2dnonconcave}
Ben{{}}{{}}{} Andrews, \emph{{Moving surfaces by non-concave curvature
  functions}}, {\tt arXiv:math.DG/0402273}.

\bibitem{AndrewsContractCalc}
Ben{{}}{} Andrews, \emph{Contraction of convex hypersurfaces in {E}uclidean
  space}, Calc. Var. Partial Differential Equations \textbf{2} (1994), no.~2,
  151--171.

\bibitem{AndrewsHarnack}
Ben{}{} Andrews, \emph{Harnack inequalities for evolving hypersurfaces}, Math.
  Z. \textbf{217} (1994), no.~2, 179--197.

\bibitem{AndrewsStones}
Ben Andrews, \emph{Gauss curvature flow: the fate of the rolling stones},
  Invent. Math. \textbf{138} (1999), no.~1, 151--161.

\bibitem{MinkowskiFlow}
Kai-Seng Chou and Xu-Jia Wang, \emph{A logarithmic {G}auss curvature flow and
  the {M}inkowski problem}, Ann. Inst. H. Poincar\'e Anal. Non Lin\'eaire
  \textbf{17} (2000), no.~6, 733--751.

\bibitem{ChowGulliverSphere}
Bennett Chow and Robert Gulliver, \emph{Aleksandrov reflection and nonlinear
  evolution equations. {I}. {T}he {$n$}-sphere and {$n$}-ball}, Calc. Var.
  Partial Differential Equations \textbf{4} (1996), no.~3, 249--264.

\bibitem{ChowLuMP}
Bennett Chow and Peng Lu, \emph{The maximum principle for systems of parabolic
  equations subject to an avoidance set}, Pacific J. Math. \textbf{214} (2004),
  no.~2, 201--222.

\bibitem{ChowAsian}
Bennett Chow and Dong-Ho Tsai, \emph{Expansion of convex hypersurfaces by
  nonhomogeneous functions of curvature}, Asian J. Math. \textbf{1} (1997),
  no.~4, 769--784.

\bibitem{ChowNonhomogGauss}
Bennett Chow and Dong-Ho{} Tsai, \emph{Nonhomogeneous {G}auss curvature flows},
  Indiana Univ. Math. J. \textbf{47} (1998), no.~3, 965--994.

\bibitem{CGFlowSpheres}
Claus Gerhardt, \emph{Flow of nonconvex hypersurfaces into spheres}, J.
  Differential Geom. \textbf{32} (1990), no.~1, 299--314.

\bibitem{CGJDG1996}
Claus{} Gerhardt, \emph{Closed {W}eingarten hypersurfaces in {R}iemannian
  manifolds}, J. Differential Geom. \textbf{43} (1996), no.~3, 612--641.

\bibitem{HamiltonThree}
Richard~S. Hamilton, \emph{Three-manifolds with positive {R}icci curvature}, J.
  Differential Geom. \textbf{17} (1982), no.~2, 255--306.

\bibitem{HuiskenRoundSphere}
Gerhard Huisken, \emph{Flow by mean curvature of convex surfaces into spheres},
  J. Differential Geom. \textbf{20} (1984), no.~1, 237--266.

\bibitem{HuiskenIlmanenPenrose}
Gerhard Huisken and Tom Ilmanen, \emph{The inverse mean curvature flow and the
  {R}iemannian {P}enrose inequality}, J. Differential Geom. \textbf{59} (2001),
  no.~3, 353--437.

\bibitem{HuiskenIlmanenIMCF2}
Gerhard Huisken and Tom{} Ilmanen, \emph{Higher regularity of the inverse mean
  curvature flow}, 2002,
  \texttt{http://www.math.ethz.ch/$\sim$ilmanen/papers/pub.html}.

\bibitem{HuiskenPolden}
Gerhard Huisken and Alexander Polden, \emph{Geometric evolution equations for
  hypersurfaces}, Calculus of variations and geometric evolution problems
  (Cetraro, 1996), Lecture Notes in Math., vol. 1713, Springer, Berlin, 1999,
  pp.~45--84.

\bibitem{IvochkinaNehringTomi}
Nina~M. Ivochkina, Thomas Nehring, and Friedrich Tomi, \emph{Evolution of
  starshaped hypersurfaces by nonhomogeneous curvature functions}, Algebra i
  Analiz \textbf{12} (2000), no.~1, 185--203.

\bibitem{JamesMcCoyArea}
James~A. McCoy, \emph{The surface area preserving mean curvature flow}, Asian
  J. Math. \textbf{7} (2003), no.~1, 7--30.

\bibitem{OSMM}
Oliver~C. Schn\"urer{}, \emph{{T}ranslating solutions to the second boundary
  value problem for curvature flows}, Manuscripta Math. \textbf{108} (2002),
  no.~3, 319--347.

\bibitem{OSA2}
Oliver~C. Schn\"urer, \emph{Surfaces contracting with speed $|{A}|^2$}, 
  J. Differential Geom. \textbf{71} (2005), no.~3, 347--363,
  {\tt arXiv:math.DG/0409388}.

\bibitem{OSKnutCrelle}
Oliver~C. Schn\"urer and Knut Smoczyk{}, \emph{{E}volution of hypersurfaces in
  central force fields}, J. Reine Angew. Math. \textbf{550} (2002), 77--95.

\bibitem{OSKnutAIHP}
Oliver~C. Schn\"urer and Knut Smoczyk, \emph{{N}eumann and second boundary
  value problems for {H}essian and {G}au{\ss} curvature flows}, Ann. Inst. H.
  Poincar\'e. Anal. Non Lin\'eaire \textbf{20} (2003), no.~6, 1043--1073.

\bibitem{OSFelixH2}
Felix Schulze{, appendix with Oliver C. Schn\"urer}, \emph{Convexity estimates
  for flows of hypersurfaces by powers of the mean curvature}, 2005, in
  preparation.

\bibitem{SmoczykAsian2000}
Knut Smoczyk, \emph{Remarks on the inverse mean curvature flow}, Asian J. Math.
  \textbf{4} (2000), no.~2, 331--335.

\bibitem{KnutHarmonicMCF}
Knut{} Smoczyk, \emph{A representation formula for the inverse harmonic mean
  curvature flow}, 2003, MPI-MIS Preprint 85/2003, {\tt
  http://www.mis.mpg.de/}.

\bibitem{Spivak4}
Michael Spivak, \emph{{A comprehensive introduction to differential geometry.
  Vol. IV. 2nd ed.}}, {Berkeley: Publish Perish, Inc. VII, 561 p.}, 1979.

\bibitem{TsoPoint}
Kaising Tso, \emph{Deforming a hypersurface by its {G}auss-{K}ronecker
  curvature}, Comm. Pure Appl. Math. \textbf{38} (1985), no.~6, 867--882.

\bibitem{UrbasExpandMZ}
John I.~E. Urbas, \emph{On the expansion of starshaped hypersurfaces by
  symmetric functions of their principal curvatures}, Math. Z. \textbf{205}
  (1990), no.~3, 355--372.

\bibitem{UrbasExpandJDG}
John I.~E.{} Urbas, \emph{An expansion of convex hypersurfaces}, J.
  Differential Geom. \textbf{33} (1991), no.~1, 91--125.

\end{thebibliography}
\def\unterstrich{\underline{\rule{1ex}{0ex}}} \def\cprime{$'$} \def\cprime{$'$}
  \def\cprime{$'$}
\providecommand{\bysame}{\leavevmode\hbox to3em{\hrulefill}\thinspace}
\providecommand{\MR}{\relax\ifhmode\unskip\space\fi MR }
\providecommand{\MRhref}[2]{%
  \href{http://www.ams.org/mathscinet-getitem?mr=#1}{#2}
}
\providecommand{\href}[2]{#2}

\end{document}